\documentstyle[11pt,psfig]{article}
\makeatletter
\def\@normalsize{\@setsize\normalsize{10pt}\xpt\@xpt
\abovedisplayskip 10pt plus2pt minus5pt\belowdisplayskip \abovedisplayskip
\abovedisplayshortskip \z@ plus3pt\belowdisplayshortskip 6pt plus3pt
minus3pt\let\@listi\@listI}
\def\subsize{\@setsize\subsize{12pt}\xipt\@xipt}
\def\section{\@startsection {section}{1}{\z@}{1.0ex plus 1ex minus
 .2ex}{.2ex plus .2ex}{\large\bf}}
\def\subsection{\@startsection {subsection}{2}{\z@}{.2ex plus 1ex}
{.2ex plus .2ex}{\subsize\bf}}
\makeatother

\newtheorem{theorem}{\sc Theorem}
\newtheorem{lemma}{\sc Lemma}
\newtheorem{coro}{\sc Corollary}
\newtheorem{nota}{\sc Notation}
\newtheorem{defin}{\sc Definition}
\newtheorem{rem}{\sc Remark}
\newtheorem{cla}{\sc Claim}
\newtheorem{ex}{\sc Example}

\newenvironment{proof}{\par \sc Proof.\rm}{\hspace*{\fill}$\Box$\vspace{1ex}}

\newenvironment{claim}{\begin{cla}}{\end{cla}}
\newenvironment{corollary}{\begin{coro}}{\end{coro}}
\newenvironment{definition}{\begin{defin}}{\end{defin}}
\newenvironment{remark}{\begin{rem}}{\end{rem}}
\newenvironment{notation}{\begin{nota}}{\end{nota}}

\title{The Average-Case Area of Heilbronn-Type Triangles\thanks{
A preliminary version of this work was presented at the 14th
Computational complexity Conference held in Atlanta in 1999.}}
\author{Tao Jiang\thanks{Supported in part by the NSERC
Research Grant OGP0046613 and a CITO grant.
Address: 
Dept of Computer Science,
University of California,
Riverside, CA 92521,
USA.
Email: jiang@cs.ucr.edu}\\
University of California, Riverside
\and
Ming Li\thanks{
Supported in part by
the NSERC Research Grant OGP0046506, a CITO grant, and the
Steacie Fellowship. Address:
Department of Computer Science,
University of California,
Santa Barbara, CA 93106,
USA.
E-mail: mli@cs.ucsb.edu}\\
University of California, Santa Barbara
\and
Paul Vit\'{a}nyi\thanks{Partially supported by the European Union
through ESPRIT BRA IV NeuroCOLT II Working Group
EP 27150.
Address: CWI, Kruislaan 413, 1098 SJ Amsterdam, The Netherlands.
Email: paulv@cwi.nl}\\
CWI and University of Amsterdam}
\date{}
 
\begin{document}
\maketitle
 
\begin{abstract}
 From among $ {n \choose 3}$ triangles with vertices chosen from 
$n$ points in the unit square, let $T$ be the one with the smallest area,
and let $A$ be the area of $T$. 
Heilbronn's triangle problem 
asks for the maximum value assumed by $A$ over all choices
of $n$ points. We consider the average-case:
If the $n$ points are chosen independently and at random (with a uniform
distribution), then
there exist positive constants $c$ and $C$ such that $c/n^3 < \mu_n < C/n^3$
for all large enough values of  $n$,
where $\mu_n$ is the expectation of $A$. Moreover,  $c/n^3 < A < C/n^3$,
with probability
close to one.
Our proof uses the incompressibility
method based on Kolmogorov complexity; it actually determines the area
of the smallest triangle for an arrangement in ``general position.'' 
\end{abstract}
 
\section{Introduction}

 From among ${n \choose 3}$ triangles with vertices 
chosen from among
$n$ points in the unit circle, let $T$ be the one of least area, 
and let $A$ be the area of $T$. Let $\Delta_n$ be the maximum assumed
by $A$ over all choices of $n$ points.
H.A. Heilbronn (1908--1975) 
asked for the exact value or approximation
of $\Delta_n$.
The list 
\cite{B99,B90,BHL97,CMS93,G72,G94,KPS81,KPS82,OPS85,R51,R72a,R72b,R73,
R76,S72,T94,YZZ91,YZZ92,YZZ94}
is a selection of papers dealing with the problem.
Obviously, 
the value of $\Delta_n$ will change only by a small
constant factor for every unit area convex shape,
and it has become customary to consider the unit square~\cite{R76}.
A brief history is as follows.
Heilbronn observed
the trivial upper bound
\footnote{We use $Omega$ in the Hardy and Littlewood sense of
``infinitely often'' as opposed to
the Knuth sense of ``always:''
If $f$ and $g$ are functions on the real numbers, then
$f(x)=O(g(x))$
if there are constants $c,x_0 > 0$ such that
$|f(x)|  \leq  c | g(x)|$, for all $x \geq x_0$;
$f(x)=o(g(x))$ if
$\lim_{x \rightarrow  \infty } f(x)/g(x) = 0$;
$f(x)= \Omega (g(x))$
if $f(x) \neq o(g(x))$.
$f(x)= \Theta (g(x))$ if both $f(x)=O(g(x))$
and $f(x)= \Omega (g(x))$.
 }
$\Delta_n = O(1/n)$ and conjectured that
$\Delta_n  =  O (1/n^2)$, and
P. Erd\H{o}s proved that this conjecture---if true---would
be tight since $\Delta_n = \Omega(1/n^2)$ \cite{R51}.
The first nontrivial result due to K.F. Roth in 1951 
established the upper bound $\Delta_n = O(1/(n\sqrt{\log \log n}))$
\cite{R51}, which was improved in
1972 by W.M. Schmidt to $O(1/(n\sqrt{\log n}))$ \cite{S72} and 
in the same year
by Roth first 
to $O(1/n^{1.105 \ldots})$~\cite{R72a} and then to
$\Delta_n = O(1/n^{1.117 \ldots})$ ~\cite{R72b}.
Roth simplified his arguments in 1973 and 1976 \cite{R73,R76}.
Exact values of $\Delta_n$ for $n \leq 15$
were studied in \cite{G72,YZZ91,YZZ92,YZZ94}.
In 1981, J. Koml\'os, J. Pintz, and E. Szemer\'edi 
\cite{KPS81}
improved Roth's upper bound to
$O(1/n^{8/7 - \epsilon })$, using the simplified arguments of Roth.
The really surprising news came in 1982 when  the same authors \cite{KPS82}
derived a lower bound $\Omega(\log n/n^2)$, narrowly refuting
Heilbronn's original conjecture.
Some believe that this lower bound is perhaps 
the best possible \cite{Er85,EP95}.
In 1997 C. Bertram-Kretzberg, T. Hofmeister, 
and H. Lefmann~\cite{BHL97} gave an algorithm that finds
a specific set of $n$ points in the unit square 
whose $\Delta_n$ (as defined above) is $\Omega(\log n/n^2)$
for every fixed $n$, using a discretization of the problem. In 1999
G. Barequet \cite{B99} derived lower bounds on  
$d$-dimensional versions of Heilbronn's problem
where $d>2$. 
All of this work concerns the {\em worst-case} value of the minimal
triangle area. 

{\bf Results:} Here we consider the {\em expected} value:
If the $n$ points are chosen independently and at random (uniform
distribution) then
there exist positive constants $c$ and $C$ such that $c/n^3 < \mu_n < C/n^3$
for all large enough $n$,
where $\mu_n$ is the expectation of the area $A$ of the smallest triangle
formed by any three points.
Moreover, with probability
close to one, $c/n^3 < A < C/n^3$.
This follows directly from corollaries~\ref{expectedlowerbound} and
\ref{expectedupperbound} of
Theorems~\ref{lowerbound} and \ref{upperbound}.
Our technique is to discretize the problem and show that all Kolmogorov-random
arrangements (see below) of $n$ points in the unit square satisfy this range
of area of the smallest triangle, where the constants $c,C$
are functions of the ``randomness deficiency'' of the arrangement---that is,
how far the Kolmogorov complexity of the arrangement falls short of
the maximum attainable Kolmogorov complexity. A Kolmogorov-random
arrangement is a rigorous way to say that the arrangement is in
``general position'' or ``typical'': there are no simple describable properties
that can distinguish any such arrangement from another one \cite{LiVi93}.
As a consequence, every arrangement in which the smallest triangle
has area outside this range---smaller or larger---cannot be Kolmogorov
random. 
According to a recent article \cite{NS99}, this result can act as
a mathematical guarantee of the afficacy of certain pseudo
Monte Carlo methods to determine the fair market value of derivatives
(on the stock market)---these methods give a sequence of points
satisfying certain pseudo-randomness properties but having less clustering
and larger smallest triangles than to be expected from truly random
sequences.
For its use in geometrical modeling see \cite{B99}.

{\bf Technique:} Our analysis uses
the {\em incompressibility method} based on Kolmogorov complexity.
The argument proceeds by using some property to be contradicted
to obtain a short encoding for some object. In the present
paper the object concerned
is usually an arrangement of $n$ pebbles on a $K \times K$ grid. 
The Kolmogorov complexity of the object is a lower bound
on the length of an encoding of the object. A contradiction arises
by the short encoding having length below the Kolmogorov complexity.
We have found that thinking in terms of coding 
is often helpful to solve our problems.
Afterwards, there may arise alternative proofs using counting,
as in the case of \cite{JLV98}, or the probabilistic method
with respect to the present result
\footnote{John Tromp has informed us in December 1999 that, 
following a preliminary version \cite{JLV99} of this work,
he has given an alternative proof of the main result
 based on the probabilistic method.
}. 
In some cases \cite{JSV97} no other proof methods seem to work.
Thinking in terms of
code length and Kolmogorov complexity enabled advances in 
problems that were open for decades, like for example \cite{JSV97,JLV98}. 
Although the technique has been
widely used in a plethora of applications, see the survey \cite{LiVi93},
it is not yet as familiar as the counting method or the probabilistic method.
One goal of the present paper
is to widen acquaintance with it
by giving yet another nontrivial example of its application. 

\section{Kolmogorov Complexity and the Incompressibility Method}

We give some definitions to establish notation.
For introduction, details, and proofs, see \cite{LiVi93}.
We write {\em string} to mean a finite binary string.
  Other finite objects can be encoded into strings in natural
ways.
  The set of strings is denoted by $\{0,1\}^*$. The {\em length}
of a string $x$ is denoted by $l(x)$, distinguishing it
from the {\em area} $|PQR|$ of a triangle on the points $P,Q,R$ in the plane.

Let $x,y,z \in \mathcal{N}$, where
$\mathcal{N}$ denotes the set of natural
numbers.
Identify
$\mathcal{N}$ and $\{0,1\}^*$ according to the
correspondence
 \[
 (0, \epsilon ), (1,0), (2,1), (3,00), (4,01), \ldots .
 \]
Here $\epsilon$ denotes the {\em empty word} with no letters.
The {\em length} $l(x)$ of $x$ is the number of bits
in the binary string $x$. 

The emphasis is on binary sequences only for convenience;
observations in any alphabet can be so encoded in a way
that is `theory neutral'.

{\bf Self-delimiting Codes:}
A binary string $y$
is a {\em proper prefix} of a binary string $x$
if we can write $x=yz$ for $z \neq \epsilon$.
 A set $\{x,y, \ldots \} \subseteq \{0,1\}^*$
is {\em prefix-free} if for any pair of distinct
elements in the set neither is a proper prefix of the other.
A prefix-free set is also called a {\em prefix code}.
Each binary string $x=x_1 x_2 \ldots x_n$ has a
special type of prefix code, called a
{\em self-delimiting code},
\[ \bar x = 1^n 0 x_1x_2 \ldots x_n  .\]
This code is self-delimiting because we can determine where the
code word $\bar x$ ends by reading it from left to right without
backing up. Using this code we define
the standard self-delimiting code for $x$ to be
$x'=\overline{l(x)}x$. It is easy to check that
$l(\bar x ) = 2 n+1$ and $l(x')=n+2 \log n+1$.

Let $\langle \cdot ,\cdot \rangle$ be a standard one-one mapping
from $\mathcal{N} \times \mathcal{N}$
to $\mathcal{N}$, for technical reasons chosen such that
$l(\langle x, y \rangle) = l(y) + l(x) + 2 l( l(x)) +1$, for example
$\langle x, y \rangle = x'y =
1^{l(l(x))}0l(x)xy$.

{\bf Kolmogorov Complexity:}
 Informally, the Kolmogorov complexity, or algorithmic entropy, $C(x)$ of a
string $x$ is the length (number of bits) of a shortest binary 
program (string) to compute
$x$ on a fixed reference universal computer 
(such as a particular universal Turing machine).
  Intuitively, $C(x)$ represents the minimal amount of information
required to generate $x$ by any effective process, \cite{Ko65}.
  The conditional Kolmogorov complexity $C(x \mid y)$ of $x$ relative to
$y$ is defined similarly as the length of a shortest program
to compute $x$, if $y$ is furnished as an auxiliary input to the
computation.
  The functions $C( \cdot)$ and $C( \cdot \mid  \cdot)$,
though defined in terms of a
particular machine model, are machine-independent up to an additive
constant (depending on the particular enumeration of Turing machines
and the particular reference universal Turing machine selected).
 They acquire an asymptotically universal and absolute character
through Church's thesis,
 and from the ability of universal machines to
simulate one another and execute any effective process,
 see for example \cite{LiVi93}.
Formally:

\begin{definition}\label{def.KC}
\rm
Let $T_0 ,T_1 , \ldots$ be a standard enumeration
of all Turing machines.
Choose a universal Turing machine
$U$ that expresses its universality in the following manner:
\[U(\langle \langle i,p \rangle ,y \rangle ) =
T_i (\langle p,y\rangle) \]
 for all $i$ and $\langle p,y\rangle$, where $p$ denotes a Turing
program for $T_i$ and $y$ an input.
We fix $U$ as our {\em reference universal computer} and define
the {\em conditional Kolmogorov complexity} of $x$ given $y$
by
\[C(x \mid y) = \min_{q \in \{0,1\}^*} \{l(q): U (\langle q,y\rangle)=x \},
 \]
for every $q$ (for example $q= \langle i,p \rangle$ above) and 
auxiliary input $y$.
The {\em unconditional Kolmogorov complexity} of $x$ is defined
by $C(x) = C(x \mid  \epsilon )$. For convenience we write 
$C(x,y)$ for $C(\langle x,y \rangle)$, and $C(x \mid y,z)$
for $C(x \mid \langle y,z \rangle)$.
\end{definition}
{\bf Incompressibility:}
Since there is a Turing machine, say $T_i$, that computes the identity
function $T_i (x) \equiv x$, it follows that
$U(\langle i,p \rangle ) =
T_i (p)$.
Hence, $C(x)  \leq l(x) + c$ for 
fixed $c \leq  2 \log i +1$ and all $x$. 
\footnote{``$2 \log i$'' and not ``$\log i$'' since we need to encode $i$
in such a way that $U$ can determine the end of the encoding. One way to do
that is to use the code $1^l(l(i))0l(i)i$ which has length $2l(l(i))+l(i)+1
< 2 \log i$ bits.
}
\footnote{In what follows, ``$\log$'' denotes the binary logarithm.
 ``$\lfloor r \rfloor$'' is the greatest integer $q$
such that $q \leq r$.
}

It is easy to see that there are also strings that can be described
by programs much shorter than themselves. For instance, the
function defined by $f(1) = 2$ and $f(i) = 2^{f(i-1)}$
for $i>1$ grows very fast, $f(k)$ is a ``stack'' of $k$ twos.
Yet for every $k$ it is clear that $f(k)$
has complexity at most $C(k) + O(1)$.
What about incompressibility?
For every $n$ there are $2^n$ binary
strings of length $n$, but only
$\sum_{i=0}^{n-1} 2^i = 2^n -1$ descriptions in binary string
format of length less than $n$.
Therefore, there is at least one binary string
$x$ of length $n$ such that $C(x)   \geq   n$.
We call such strings $incompressible$. The same argument holds
for conditional complexity: since for every length $n$
there are at most $2^n-1$ binary programs of length $<n$,
for every binary string $y$
there is a binary string $x$ of length $n$ such that
$C(x \mid  y)   \geq   n$. 
Strings that are incompressible 
are patternless,
since a pattern could be used to reduce
the description length. Intuitively, we
think of
such patternless sequences
as being random, and we
use ``random sequence''
synonymously with ``incompressible sequence.''
Since there are few short programs, there can
be only few objects of low complexity:
the number
of strings of length $n$ that are compressible by at most $\delta$ bits
is at least $2^n - 2^{n-\delta} +1$. 
\begin{lemma}
\label{C2}
Let $\delta$ be a positive integer.
For every fixed $y$, every
set $S$ of cardinality $m$ has at least $m(1 - 2^{-\delta} ) + 1$
elements $x$ with $C(x \mid  y)   \geq   \lfloor \log m \rfloor  - \delta$.
\end{lemma}
\begin{proof}
There are $N=\sum_{i=0}^{n-1} 2^i = 2^n -1$ binary strings 
of length less than $n$. A fortiori there are at most $N$
elements of $S$ that can be computed by
binary programs of length less than $n$, given $y$.
This implies that at least $m-N$  elements of $S$ cannot
be computed by binary programs of length less than $n$, given $y$.
Substituting $n$ by  $\lfloor \log m \rfloor  - \delta$ together
with Definition~\ref{def.KC} yields the lemma. 
\end{proof}

If we are given $S$ as an explicit table
then we can simply enumerate its elements
(in, say, lexicographical order) using a fixed program not depending
on $S$ or $y$. Such a fixed program can be given in $O(1)$ bits.
Hence the complexity satisfies $C(x \mid S,y) \leq \log |S| + O(1)$.

{\bf Incompressibility Method:}
In a typical proof using the incompressibility method,
one first chooses an incompressible object from the
class under discussion.
The argument invariably says that if a desired property
does not hold, then in contrast with the assumption, the object
can be compressed. This yields the required contradiction.
Since most objects are almost incompressible, the desired property
usually also holds for almost all objects, and hence on average.

\section{Grid and Pebbles}
In the analysis of the triangle problem we first
consider a discrete version based on an equally spaced $K \times K$ grid
in the unit square. The general result for the
continuous situation is then obtained by taking the limit
for $K \rightarrow \infty$. 
Call the  resulting axis-parallel $2K$ lines {\em grid lines} and
their crossing points {\em grid points}. 
We place $n$ points on grid points. These $n$ points will be referred to
as {\em pebbles} to avoid confusion with grid points or
other geometric points
arising in the discussion. 

There are ${K^2 \choose n}$ ways to put $n$ 
{\em unlabeled} pebbles on the grid where
at most one pebble is put on every grid point.
We count only distinguishable arrangements
without regard for the identities of the placed pebbles. 
Clearly, the restriction 
that no two pebbles
can be placed on the same grid point is no restriction anymore when we let $K$
grow unboundedly.

Erd\H{o}s \cite{R51}
demonstrated that for the special case of $p \times p$ grids,
where $p$ is a prime number, there are necessarily arrangements
of $p$ pebbles with every pebble placed on a grid point
such that no three pebbles are collinear. The least
area of a triangle in such an arrangement is at least $1/(2p^2)$. This
implies that the triangle constant $\Delta_n = \Omega(1/n^2)$ as
$n \rightarrow \infty$ through the special sequence of primes. 

We now give some detailed examples---used later---of
the use of the incompressibility method.
By Lemma~\ref{C2}, for every integer $\delta$ independent of $K$,
every arrangement $X_1, \ldots , X_n$ (locations of pebbles), out of at least
a fraction of $1-1/2^{\delta}$ of all arrangements of $n$ pebbles 
on the grid, satisfies
\begin{equation}\label{eq.inc1}
C(X_1, ..., X_n \mid n,K ) \geq \log {K^2 \choose n} - \delta .
\end{equation}

\begin{notation}
\rm
For convenience we abbrieviate the many occurrences of the phrase
``Let $X_1, \ldots , X_n$ be an arrangement of $n$ pebbles on the $K \times K$
grid, let $n$ be fixed and $K$ be sufficiently large,
and let $\delta$ be a positive integer constant such that
(\ref{eq.inc1}) holds''
to ``If (\ref{eq.inc1}) holds'' in the remainder of the paper.
\end{notation}

Note that, for every arrangement $X_1, \ldots , X_n$ of $n$ pebbles on a
$K \times K$ grid, we have 
$C(X_1, ..., X_n \mid n,K ) \leq \log {K^2 \choose n}+O(1)$---
there is a fixed program of $O(1)$ bits
for the reference universal computer
that reconstructs the $X_1, \ldots , X_n$ from $n,K$ and its index in
the lexicographical ordering of all possible arrangements.  
That (\ref{eq.inc1}) holds with $\delta$ small means that the arrangement
$X_1, \ldots , X_n$ of pebbles on the grid has no regularity that can
be used to prepare a description that is significantly shorter than
simply giving the index in the lexicographical ordering of all
possible choices of $n$ positions from the available $K \times K$
grid positions. We can view such an arrangement as being ``random''
or ``in general position.'' 

\begin{lemma}\label{cl.noline}
If (\ref{eq.inc1}) holds, then
no three pebbles can be collinear, and
so the area of a smallest triangle is at least $1/(2(K-1)^2)$.
\end{lemma}

\begin{remark}
\rm
This is the first proof of the paper using the incompressibility argument.
Let us explain the proof idea in detail: 
On the one hand, we construct a description $d$ such that
the arrangement $X_1, \ldots , X_n$
can be reconstructed from $d$ by a fixed program $p$ for the 
universal reference computer, given also $n$ and $K$.
If $p$ is in self-delimiting format, then the universal reference computer
can parse $pd$ into its constituent parts $p$ and $d$, and then
execute $p$ to reconstruct $X_1, \ldots , X_n$ from the auxiliary 
information $n,K$,  together with the description $d$.
On the other hand, by definition the Kolmogorov complexity of an object is the
length of its {\em shortest} program for the reference universal computer
and we have assumed a lower bound on the Kolmogorov complexity.
Since the description $pd$ is a program for the reference universal computer,
its length $l(pd)$ must be at least as large as the 
Kolmogorov complexity (the auxiliary information $n,K$ being the same
in both cases).
By the lower bound (\ref{eq.inc1}) this shows that 
$l(pd) \geq \log {K^2 \choose n} - \delta$. Since $l(p)$ is independent
of $n,K$ we can set  $l(p)=O(1)$ in this context,
and obtain $l(d) \geq \log {K^2 \choose n} - \delta -O(1)$.  
By exploiting collinearity of pebbles in the description $d$,
to make it as compact as possible, this inequality will yield
the required contradiction for $n$ fixed and $K$ large enough.

\end{remark}

\begin{proof} 
Place $n-1$ pebbles at positions chosen from the total of $K^2$ grid
points---there are ${{K^2} \choose {n-1}}$ choices. Choose two pebbles,
$P$ and $Q$, from among the $n-1$ pebbles---there are
${{n-1} \choose 2}$ choices. Choose a new pebble $R$ on the straight
line determined by $P,Q$.
The number of grid points 
on this line between $P$ (or $Q$) and $R$,
which number is $<K$, identifies $R$ uniquely 
in $ \leq \log K$ bits.
There is a fixed algorithm that, on input  $n$ and $K$,
decodes a binary
description consisting of the items above---each encoded as the logarithm
of the number of choices---and
computes the positions of the $n$
pebbles.
By (\ref{eq.inc1}) this implies 
\[  \log {{K^2} \choose {n-1}} + \log {{n-1} \choose 2}  +  \log K
+ O(1) 
  \geq \log {{K^2} \choose n} - \delta  .\]
Using the asymptotic expression
\begin{equation}\label{eq.ass}
\log  {a \choose b} 
-  b \log \frac{a}{b} \rightarrow b \log e - \frac{1}{2} \log b + O(1)
\end{equation}
for $b$ fixed and $a \rightarrow \infty$,
one obtains $3 \log n \geq \log K - \delta  + O(1)$,
which is a contradiction for $n$ fixed and $K$ sufficiently large.
\end{proof}

\begin{lemma}\label{cl.nogridline}
If (\ref{eq.inc1}) holds,
then no two pebbles can be on the same (horizontal) grid line.
\end{lemma}

\begin{proof}
Place $n-1$ pebbles at positions chosen from the total of $K^2$ grid
points---there are ${{K^2} \choose {n-1}}$ choices. Choose one pebble
$P$ from among the $n-1$ pebbles---there are
$n-1$ choices. Choose a new pebble $R$ on the (horizontal) grid line
determined by $P$---there are $K-1$ choices. There is a fixed
algorithm that, on input $n$ and $K$, 
reconstructs the positions of all $n$ pebbles
from a description of these choices.
By (\ref{eq.inc1}) this implies
\[ \log {{K^2} \choose {n-1}}+ \log (n-1)  +  \log K + O(1) 
\geq \log {{K^2} \choose n} - \delta .\]
Using (\ref{eq.ass}) with fixed $n$ and $K \rightarrow \infty$  we obtain
$ 2 \log n \geq \log K - \delta + O(1)$,
which is a contradiction for
large enough $K$.
\end{proof}


\section{Lower Bound}
Our strategy is to show that if we place $n$ pebbles
on a $K \times K$ grid, such that the arrangement has high 
Kolmogorov complexity, then every three pebbles form
a triangle of at least a certain size area. If the area
is smaller, then this can be used to compress the description size
of the arrangement to below the assumed Kolmogorov complexity.

\begin{theorem}\label{lowerbound}
If (\ref{eq.inc1}) holds,
then there is a positive constant $c_1$ such that
the least area of every triangle formed by three pebbles on the grid
is at least
$c_1/(2^{\delta}n^3)$.
\end{theorem}

\begin{proof}
Place $n-1$ pebbles at positions chosen from the total of $K^2$ grid 
points---there are ${{K^2} \choose {n-1}}$ choices. Choose two pebbles,
$P$ and $Q$, from among the $n$ pebbles---there are 
${{n} \choose 2}$ choices. Place a new pebble $R$ at one of the remaining
grid points. Without loss of generality, let
the triangle $PQR$ have $PQ$ as the longest side.
Center the grid coordinates on $P=(0,0)$ with $Q=(q_1,q_2)$
and $R=(r_1,r_2)$ in units of $1/(K-1)$ in both axes directions.
Then $R$ is one of the grid points on the
two parallel line segments of length $L=|PQ|= \sqrt{q_1^2 + q_2^2}/(K-1)$ at
distance $H = |q_2r_1-q_1r_2|/((K-1) \sqrt{q_1^2 + q_2^2})$ 
from the line segment $PQ$,
as in
Figure~\ref{lowerboundfig}.
The number of grid points on each of these 
line segments (including one endpoint and excluding the other endpoint) 
is a positive integer
$g = \gcd(q_1,q_2)$---the line $q_2 x = q_1 y $ 
has $g$ integer coordinate points between 
$(0,0)$ and $(q_1,q_2)$ including one of the 
endpoints. This implies that $f$ defined by $LH (K-1)^2 = fg$ is
a positive integer as well.

\begin{figure}[htbp]
\hfill\ \psfig{figure=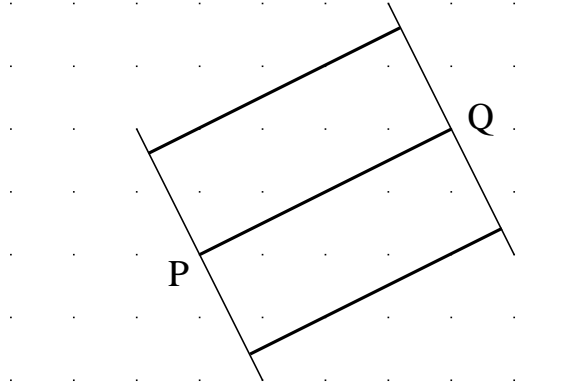,width=3.0in} \hfill\
\caption{Smallest triangle based on pebbles $P,Q$.}
\label{lowerboundfig}
\end{figure}

Enumerating the grid points concerned in lexicographical order, 
the index of $R$ takes at most 
$\log (2gf) = \log (2g) + \log f= \log (4|PQR| (K-1)^2)$
bits, where $|PQR|$ denotes the area of the triangle PQR. 
Altogether this constitutes an effective description of the arrangement
of the $n$ pebbles.
By the assumption in the theorem the arrangement satisfies (\ref{eq.inc1}),
that is, the number of bits
involved in any effective description of the arrangement
is lower bounded by the righthand side.
Then,
\[ \log {{K^2} \choose {n-1}} + \log {{n} \choose 2} +
 \log ( 4 |PQR| (K-1)^2) +  O(1) \geq
\log {K^2 \choose n} - \delta . 
\] 
By approximation (\ref{eq.ass}),
\[ \log {K^2 \choose n} - \log {{K^2} \choose {n-1}} \rightarrow
\log \frac{K^2}{n} + O(1) \]
for large enough fixed $n$ and $K \rightarrow \infty$.
Therefore,
 $ \log |PQR| + O(1)  \geq - 3 \log n - \delta  + O(1), K \rightarrow 
\infty$.
Consequently, there exists a positive constant $c_1$, independently
of the particular triangle $PQR$, such that $|PQR| > c_1/(n^3 2^{\delta })$
for all large enough $n$ and $K$. Since this holds for every triangle
$PQR$, constructed as above,
 it holds in particular for a triangle of least area $A$.
\end{proof}

By Lemma~\ref{C2}
the probability concentrated on the set of arrangements satisfying
(\ref{eq.inc1}) is at least $1- 1/2^{\delta }$:
\begin{corollary}
If $n$ points are chosen independently and at random (uniform distribution) 
in the unit square, and $A$ is the least area of a triangle formed by three
points, then there is a positive constant $c_1$ such that for every
positive $\delta$ we have $A > c_1/(2^{\delta}n^3)$
with probability at least $1-1/2^{\delta}$.
\end{corollary}
In the particular case of $\delta = 1$
the probability concentrated on arrangements satisfying
(\ref{eq.inc1}) is at least $\frac{1}{2}$ which immediately implies:
\begin{corollary}\label{expectedlowerbound}
If $n$ points are chosen independently and at random (uniform distribution) 
in the unit square, then there is a positive constant $c$ such that
the least area of some triangle formed by three points
has expectation $\mu_n > c/ n^3$.
\end{corollary}





\section{Upper Bound}
Every pair of pebbles out of an incompressible
arrangement of $n$ pebbles on a $K \times K$ grid
defines a distinct line by Lemma~\ref{cl.noline}.
The two pebbles defining such a line together with any other pebble forms
a triangle.
If $A$ is the least area of a triangle formed by three pebbles, then
this constrains the possibilities of placing a third pebble close
to a line defined by two pebbles.
Thus, every such line 
defines a forbidden strip 
on both sides of the line where no pebbles can be placed.
It is easy geometry to see that every forbidden strip
covers an interval of length $2A$ of every grid line on both sides of the
intercept of the ``forbidding'' line concerned.
Our strategy is as follows: Divide the pebbled unit square by a 
straight line parallel to the horizontal sides into two
parts containing about one half of the $n$ pebbles each.
Show that the pebbles in the larger half (the halves may not have equal
area) of the unit square 
define $\Omega (n^2)$ distinct  ``forbidding lines'', 
that cross both the dividing line and
the opposite parallel side of the unit square. 
While the associated forbidden grid point positions can overlap, 
we show that they
don't overlap too much. As a consequence the
set of grid points allowed
to place
the remaining $n/2$ pebbles in the smaller remaining half of
the unit square, gets restricted
to the point that the description of the arrangement can be
compressed too far.
This argument
is so precise that for small 
$\delta$ in (\ref{eq.inc1})
the upper bound is the same order of magnitude as the
previously proven lower
bound.

\begin{theorem}\label{upperbound}
If (\ref{eq.inc1}) holds with
$\delta  < (2- \epsilon ) \log n$ for some positive constant $\epsilon$, then
there is a positive constant $C_1$ such that the least area of some
triangle formed by three pebbles on the grid
is at most
\begin{equation}\label{eq.Udelta}
A(\delta)= \frac{14 \delta  + O(1)}{4C_1n^3 \log e}.
\end{equation}
\end{theorem}

\begin{proof}
Choose $n$ pebbles at positions chosen from the total of $K^2$
grid points such that (\ref{eq.inc1}) is satisfied.
Divide the unit square by a horizontal grid line
into an upper and a lower half, each of which
contains $n/2 \pm 1$ pebbles---there
are no grid lines containing two pebbles by Lemma~\ref{cl.nogridline}.
We write {\em forbidding line} for
a line determined by two pebbles in the
upper half that intersects
all horizontal grid lines in the lower half of the unit square.

\begin{claim}\label{cl.n2lines}
\rm
If (\ref{eq.inc1}) holds, then there is a positive 
constant $C_1$ such that there are
at least $C_1 n^2$ forbidding lines.
\end{claim}

\begin{proof}
Take the top half
to be the larger half so that it has area at least $1/2$. Divide the top
half into five vertical strips of equal width of $1/5$ and 
five horizontal strips of equal width
$1/10$ starting from the top---ignore the possibly remaining
horizontal strip at the bottom of the top half. 
Clearly, a forbidding line determined 
by a pebble in the upper rectangle and a pebble in the lower
rectangle of the middle vertical strip intersects 
the bottom horizontal grid line. We show that these rectangles
contain at least $n/100$ points each, and hence the claim
holds with $C_1 = 1/10,000$.

Consider either rectangle (the same argument will hold for
the other rectangle).
Let it contain $m \leq n$ pebbles. 
Since the area of the rectangle
is $ 1/5 \times 1/10 $ it contains $K^2/50$ 
grid points (plus or minus the grid points
 on the circumference of length $3K/5$ which we ignore).
Place $n-m$ pebbles at positions chosen from $49K^2/50$ grid
points outside the rectangle---there 
are ${{49K^2/50} \choose {n-m}}$ choices---and
place $m$ pebbles at positions chosen from the total of $K^2/50$ grid
points in the rectangle---there are ${{K^2/50} \choose {m}}$ choices. 
Given $n$ and $K$, the $n$ pebble positions are determined by $m$, 
the position
of the rectangle and an index number $i$ of $\log i$ bits with
\begin{eqnarray*}
\log i & = & \log {{49K^2/50} \choose {n-m}}
 {{K^2/50} \choose m} \\
& \rightarrow &
(n-m) \log \frac{49K^2/50}{n-m}  + m \log \frac{K^2/50}{m}
+ n \log e -  \frac{1}{2}\log nm + O(1), 
\end{eqnarray*}
for $K \rightarrow \infty$ with $n,m$ fixed,
by (\ref{eq.ass}).
Given $n$ we can describe $m$ in $\log n$ bits. Thus,
given $n$ and $K$, the total description 
length of the description of the arrangement
of the $n$ pebbles is $\log n+\log i+O(1)$ bits. This must
be at least the Kolmogorov complexity of the arrangement. Then,
by (\ref{eq.inc1}),
\[
(n-m) \log \frac{49K^2/50}{n-m} + m \log \frac{K^2/50}{m} - \frac{1}{2} \log m
+ O(1) \geq n \log \frac{K^2}{n} - \delta . \]
This implies
\begin{eqnarray*}
 \delta & \geq & (n-m) \log \frac{50(n-m)}{49} + m \log 50m + \frac{1}{2}\log m
- n \log n - O(1) \\
& > & (n-m) \log \frac{50(n-m)}{49} + m \log 50m - n \log n - O(1).
\end{eqnarray*}
Assume, by way of contradiction, $m \leq n/100$.  Then,
\begin{eqnarray*}
\delta & \geq & \frac{99}{100} n \log \frac{4950}{4900} n + 
\frac{1}{100}n \log \frac{50}{100} n - 
n \log n  - O(1) \\
& = & n ( \frac{99}{100} \log \frac{4950}{4900} 
+ \frac{1}{100} \log \frac{50}{100}) - O(1) \\
& > & n ( 0.0145 - 0.01) - O(1),
\end{eqnarray*}
which contradicts
$\delta = O( \log n)$ in the statement of the theorem.
Hence the top rectangle and the bottom
rectangle of the middle strip in the top half contain
at least $n/100$ pebbles each. Each pair of pebbles, one in the top rectangle
and one in the bottom rectangle, determine a distinct forbidding line by
Lemma~\ref{cl.noline} (no three pebbles can be collinear under
assumption (\ref{eq.inc1})).
The claim is proven with $C_1 = (1/100) \cdot (1/100) = 1/10^4$.
\end{proof}

\begin{claim}\label{cl.uppline}
\rm
Let $w_1,w_2,w_3,w_4,w_5$ be the spacings between the six
consecutive intercepts of a sextuplet of forbidding lines
with a horizontal grid line in the bottom half containing a pebble,
and let $D=w_1+w_2+w_3+w_4+w_5$.
If (\ref{eq.inc1}) holds,
then there is a positive $C_2$ such that 
 $D > C_2/n^{3- \epsilon/5}$ with $\epsilon$ as in the statement
of the theorem.
\end{claim}

\begin{proof}
Place $n-5$ pebbles at positions chosen from the total of $K^2$ grid
points---there are ${{K^2} \choose {n-5}}$ choices. Choose eight pebbles,
$P_i$ ($i=0,1,2,3,5,7,9,11$) 
from among the $n-5$ pebbles---there are
at most ${{n-5} \choose 8}$ choices---and five new pebbles 
$P_j$ ($j=4,6,8,10,12$) such that
$P_1P_2, P_3P_4, P_5P_6$, 
$P_7P_8, P_9P_{10}$,
$P_{11}P_{12}$ is the sextuplet of forbidding lines in the claim, and
$P_0$ is a pebble in the lower half.
Without loss of generality we assume that the ``middle'' pebbles 
of unknown position $P_j$ ($j=4,6,8,10,12$), as well as $P_2$
in known position, are in between
the other defining pebble of the forbidding line concerned and its
intercept with the lower grid line containing $P_0$. That is,
the top-to-bottom order a forbidding line is 
$P_1,P_2,\mbox{\rm intercept}_1$,
$P_3,P_4,\mbox{\rm intercept}_2$, and so on. Then,
a forbidding line determined by an outermost pebble and an intercept,
together with the grid line containing the middle pebble, enables us
to determine the grid point on which the middle pebble is located.
An error in the position of the intercept leads to a smaller error
in the position of the middle pebble. Thus, a precision of the position
of the intercept up to $1/(4(K-1))$, together with the precise
position of the outermost pebble, enables us to determine the
grid point containing the middle pebble as 
the unique grid point in a circle with radius
$1/(4(K-1))$ centered on the computed geometric point.   
The coordinates of the five unknown $P_j$'s
are determined by (i) the locations
of the five intercepts of the associated quintuplet of forbidding lines 
with the lower half horizontal grid
line on which $P_0$ is located, 
and (ii) the five unknown distances between these
intercepts and the $P_j$'s along the five associated forbidding lines.
The grid point positions of the $P_j$'s are uniquely determined if we know
the latter distances up to precision $1/4(K-1)$.  All six intercepts
in the statement of the theorem are in an interval
of length $D$ which contains $DK$ grid points (rounded to the
appropriate close entier value). We can describe every intercept
in this interval (up to the required precision) in $\log DK + O(1)$
bits. 
 Relative to the intersection
of the known forbidding line $P_1P_2$, therefore, item (i) uses  
$5 \log DK + O(1)$ bits. Item (ii) uses $5 \log K + O(1)$ bits. 
Given $n,K$, we can describe the placement of the $n-5$ pebbles
in $\log {{K^2} \choose {n-5}}$ bits; the choice of the eight pebbles 
among them in $\log {{n-5} \choose 8}$ bits; and we have shown
that the placement of
the five unknown pebbles can be reconstructed from an additional
$5 \log DK  +5 \log K + O(1)$ bits.  
Together this forms a description of the complete arrangement.
By (\ref{eq.inc1}) this implies:
\[ 
 \log {{K^2} \choose {n-5}} +8\log n +  5\log DK + 5\log K + O(1) \geq
\log {{K^2} \choose n} - \delta  .\] 
A now familiar calculation using (\ref{eq.ass}) yields
$  5 \log D + O(1) \geq - 13 \log n  - \delta$, for fixed $n$ and 
$ K \rightarrow \infty$.
This shows
$D > C_2 2^{( 2 \log n - \delta )/5} /n^3 $ for some positive constant $C_2$.
Substituting
$\delta < (2-\epsilon)\log n$
proves the claim. 
\end{proof}

We have now established that there are $C_1 n^2$ distinct
forbidding lines (with $C_1$ as in Claim~\ref{cl.n2lines})
determined by pairs of pebbles in the upper half,
and by construction every such forbidding line intersects every
lower half horizontal grid line. Moreover,
every $D$-length interval (with $D$ as in Claim~\ref{cl.uppline}) on 
a lower half horizontal grid line---that contains a pebble---contains
at most six intercepts of forbidding lines.  
This means that we can select $C_1 n^2 /7$ consecutive intercepts
on such a grid line
that are separated by intervals of at least length $D$. The two pebbles
$P,Q$ defining the forbidding line $l_1$, 
together with any pebble $R$ on a lower
half horizontal grid line $l_2$, determine a triangle. 
If $d$ is the distance between the intercept point of $l_1$ with $l_2$ and
the pebble $R$, and $\alpha$ is the angle between the forbidding line $l_1$
and grid line $l_2$, then the triangle side located on the forbidding line 
has length
$\leq 1/\cos \alpha$ while the height of the triangle with respect
to that side is $d \cos \alpha$.    
Thus, if $A$ is the area of the smallest triangle formed by any
three pebbles, then $d \geq 2A$.
Consequently, all grid positions in intervals of length $2A$ 
on both sides of an intercept
of a forbidding line with a lower half grid line---that contains a pebble---
are forbidden for pebble placement. As long as $2d \leq D$, or
\begin{equation}\label{eq.f2}
4A \leq D, 
\end{equation}
this means that the $C_1 n^2 /7$ consecutive intercepts
exclude $4A C_1 n^2 /7$ grid positions from pebble placement
on the horizontal lower grid line concerned.
If (\ref{eq.f2}) does not hold, that is, $4A > D$, then
at least $DC_1 n^2 /7$ grid positions are excluded.
Given the pebbles in the upper half, and therefore the forbidding lines,
the excluded grid points in the lower half are determined.
Therefore, with
\begin{equation}\label{eq.B}
 B = \min \{4A, D\} 
\end{equation}
and also given the horizontal lower half grid line concerned,
we can place a
pebble on the grid line in at most
\begin{equation}\label{eq.excl}
         K(1 - C_1 n^2B/7)
\end{equation}
positions.
We now use this fact to construct a short encoding
of the total arrangement of the $n$ pebbles satisfying (\ref{eq.inc1}):
Select $n$ horizontal grid lines (there can be only one pebble per
grid line by Lemma~\ref{cl.noline}) chosen from the total of $K$ grid
lines---there are ${{K} \choose {n}}$ choices. Select on everyone of the
upper $n/2$ horizontal grid lines a grid point to place a pebble---there
are $K^{n/2}$ choices. Finally, select in order from top to bottom on
the lower $n/2$ horizontal grid lines 
$n/2$ grid points to place the pebbles---there
are only $(K(1 - C_1 n^2B/7))^{n/2}$ choices by (\ref{eq.excl}).
Together these choices form a description of the arrangement.
Given the values of $n,K$ we can encode these choices in self-delimiting
items, and 
by (\ref{eq.inc1}) this implies:
\[
\log {K \choose n} 
 +  
\frac{n}{2} \log K + \frac{n}{2} \log K(1 - C_1n^2B/7) + O(1) \geq
\log {K^2 \choose n} - \delta .
\]
Using (\ref{eq.ass})
with $n$ fixed yields
\[ \frac{n}{2} \log (1 - C_1n^2B/7) \geq - \delta  - O(1), 
\; \; K \rightarrow \infty .\]
The left-hand side 
\[ \log \left( 1 - \frac{C_1n^3B/14}{n/2} \right)^{n/2}
= \log e^{-C_1n^3B/14}, \;\; n \rightarrow \infty , \]
so that
\begin{equation}\label{eq.fn}
B \leq \frac{14 \delta  + O(1)}{C_1n^3 \log e} 
\end{equation}
Since $\delta < 2 \log n$ in the right-hand side, Claim~\ref{cl.uppline} 
shows that $D>B$. Therefore,
 (\ref{eq.B}) implies $B=4A$ so that (\ref{eq.fn}) establishes the theorem.
\end{proof} 

Together with Lemma~\ref{C2}, Theorem~\ref{upperbound} 
implies that the smallest triangle
in an arrangement has an area below a particular upper bound
with a certain probability.
\begin{corollary}
If $n$ points are chosen independently and at random (uniform distribution)
in the unit square, and $A$ is the least area of a triangle formed by three
points, then for every
positive $\delta < (2-\epsilon ) \log n$ ($\epsilon > 0$),
we have 
\[ A < A(\delta)\]
with probability at least $1-1/2^{\delta}$.
\end{corollary}

That is, the probability that $A <  A(1)$
at least $\frac{1}{2}$ ($\delta = 1$),
the probability that $A < A(2)$ is at least $\frac{3}{4}$
($\delta = 2$),
and so on. Since $A(\delta+1) \geq A(\delta)$,
we can upper bound the expectation
$\mu_n$ of $A$ by upper bounding the probability of $A$ with
$A(\delta) < A \leq A(\delta+1)$ by 
$2^{-\delta+1} = [(1 - 2^{-\delta })-(1-2^{-\delta-1})]$. 
We do this for $\delta \leq  1.9 \log n $. 
The remaining probability is 
$1/n^{1.9}$ or slightly less (because $\delta$ is integer). This probability
is so small that, even if we assume the known worst-case
upper bound on $A$ for the remaining cases,
known to be $C_3/n^{8/7 - \epsilon' }$ for some positive constant $C_3$ for 
every $\epsilon ' > 0$,
\cite{KPS81}, the result is insignificant. There is a positive constant
$C$ such that:
\begin{eqnarray*}
\mu_n &\leq & 
  \sum_{\delta =1}^{ 1.9 \log n }
2^{-\delta} 
A(\delta) +
 \frac{1}{n^{1.9}} \frac{C_3}
{n^{8/7 - \epsilon ' }} < \frac{C}{n^3} .
\end{eqnarray*}
\begin{corollary}\label{expectedupperbound}
If $n$ points are chosen independently and at random (uniform distribution)
in the unit square, then there is a positive constant $C$ such that
the least area of some triangle formed by three points
has expectation $\mu_n < C/ n^3$.
\end{corollary}

\noindent 
{\bf Acknowledgement:}
We thank John Tromp for  
help with the
proof of Theorem~\ref{lowerbound} and him and the anonymous referees for
valuable comments on drafts of this paper.

\end{document}